\documentclass[11pt]{article}

\RequirePackage{amsthm}
\RequirePackage{amssymb}
\usepackage[all,ps]{xy}

\usepackage{latexsym}
\usepackage[frenchb]{babel}

\usepackage[latin1]{inputenc}

\usepackage{a4wide}

\author{Andrea Surroca Ortiz}
\title{\Large Valeurs alg\'ebriques de fonctions transcendantes}
\date{3 ao\^{u}t 2005}

\newtheorem{thm}{Th\'eor\`eme}[section]

\newtheorem{lemme}[thm]{Lemme}

\newtheorem{remarque}[thm]{Remarque}

\newcommand\rat{\mathbf{Q}}

\newcommand\rationnels{\mathbf{Q}}

\newcommand\Qbarre{\overline{\mathbf{Q}} }
  
 \newcommand\naturels{\mathbf{N}}

\newcommand\enteros{\mathbf{Z}}

\newcommand\complexes{{\mathbf{C}}}

\newcommand\DRferme{\overline{D(0,R)}}

\newcommand\Drferme{\overline{D(0,r)}}

\newcommand\Duferme{\overline{D(0,1)}}

\newcommand\DRouvert{D(0,R)}

\newcommand\Duouvert{D(0,1)}

\newcommand\Drhoferme{\overline{D(0,\rho)}}

\newcommand\card{\mathrm{card}}

\newcommand\pgcd{\mathrm{pgcd}}

\newcommand\esp{\hspace{0,2cm}}
\newcommand\espa{\hspace{0,17cm}}

\begin{document}

\maketitle

\begin{quote}
\textbf{Abstract.}
{\small We study the set of algebraic numbers of bounded height and bounded degree where an analytic transcendental function takes algebraic values. 
 }
\end{quote}

\begin{quote}
\textbf{R\'esum\'e.}
{\small  On \'etudie l'ensemble des nombres alg\'ebriques de hauteur
  et de degr\'e born\'es o\`u une fonction analytique transcendante prend des
  valeurs alg\'ebriques. }
\end{quote}


\section{ Introduction.}

\'Etant donn\'ee une fonction $f$ analytique, on consid\`ere l'ensemble $S_f$ des points alg\'ebriques en lesquels la fonction $f$ prend des valeurs alg\'ebriques. Ici, la fonction $f$ sera transcendante sur $\mathbf{C}(z)$, c'est-\`a-dire qu'il n'existe pas de polyn\^ome non nul en 2 variables et \`a coefficients complexes, s'annulant sur tous les points $(z,f(z))$ pour les $z$ o\`u $f$ est d\'efinie.  

Par exemple, d'apr\`es le th\'eor\`eme de Hermite-Lindemann (th\'eor\`eme 1.2 \cite{w}), la fonction exponentielle prend des valeurs transcendantes en tout point alg\'ebrique, sauf en $0.$ Donc pour $f(z) = e^z,$ on a $S_f = \{0\}.$ Par la m\^eme raison, si on pose $f(z) = e^{P(z)}$ o\`u $P \in \Qbarre[X],$ alors $S_f$ est l'ensemble des z\'eros de $P.$ En supposant vraie la conjecture de Schanuel (\cite{lang} p.30 ou \cite{w}), pour $f(z)=\sin(\pi z)e^z$ on a $S_f=\enteros.$ Le th\'eor\`eme de Gel'fond-Schneider (th\'eor\`eme 1.4 \cite{w}) nous fournit d'autres exemples: si $f(z) = e^{\lambda z}$ o\`u $\lambda \ne 0$ est tel que $e^{\lambda}$ soit alg\'ebrique (par exemple, $f(z) = 2^z$ ou $f(z) = e^{i\pi z}$), alors $S_f = \rationnels.$

Dans une lettre de 1886 adress\'ee \`a Strauss \cite{s1}, Weierstrass sugg\'erait l'existence de fonctions enti\`eres transcendantes prenant des valeurs alg\'ebriques en tous les points alg\'ebriques. 
%
Apr\`es un premier r\'esultat de Strauss allant dans cette direction, St\"ackel (cf. \cite{gramain} et \cite{s1}),  \'enonce le th\'eor\`eme suivant:

\emph{\'Etant donn\'es un sous-ensemble $\Sigma$ d\'enombrable de $\complexes$ et un sous-ensemble $T$ dense de $\complexes,$ il existe une fonction enti\`ere transcendante envoyant $\Sigma$ dans $T.$}

En prenant $\Sigma=T=\Qbarre,$ on obtient une fonction $f$ enti\`ere et transcendante pour laquelle l'ensemble $S_f$ est $\Qbarre$ tout entier. Par ailleurs, on en d\'eduit, en posant $\Sigma= \Qbarre$ et $T= \complexes \setminus \Qbarre,$ l'existence d'une fonction enti\`ere transcendante prenant des valeurs transcendantes en tous les points alg\'ebriques, c'est-\`a-dire, telle que $S_f= \emptyset$.  (Sous la conjecture de Schanuel, on peut montrer que pour la fonction $f(z) = e^{e^z}$ on a $S_f=\emptyset$.)
 
 F. Gramain remarque dans \cite{gramain}, que si $\Sigma$ est contenu dans $\mathbf{R}$, la m\^eme d\'emonstration s'applique avec $T$ dense dans $\mathbf{R}$. Ceci nous permet de prendre $\Sigma = K$, corps de nombres r\'eel, ou $\Sigma = O_{K,S}$ l'anneau des $S$-entiers d'un corps de nombres $K$ r\'eel, et, dans les deux cas, $T = \mathbf{Z}\left[ \frac{1}{n} \right]$, pour n'importe quel entier naturel $n\geq 2$.

En suivant une remarque de P. St\"ackel, qui construit \cite{s2} une fonction alg\'ebri\-que qui prend, ainsi que toutes ses d\'eriv\'ees, des valeurs alg\'ebriques en tous les points alg\'ebriques, G. Faber \cite{faber} construit une fonction $G$ enti\`ere et transcendante telle que, pour tout nombre entier $k$ et tout $\alpha \in \Qbarre,$ la valeur de la fonction d\'eriv\'ee $G^{(k)}$ en $\alpha$ est dans $\rat + i \rat.$

Dans la section \ref{exemples}, nous construisons une fonction $f$ enti\`ere et transcendante (th\'eor\`eme \ref{contrex}), dont toutes ses d\'eriv\'ees, envoient tout corps de nombres dans lui-m\^eme (cf. aussi \cite{vdp}), et pour laquelle on contr\^ole la hauteur des valeurs prises aux points alg\'ebriques.
On construit aussi (th\'eor\`eme \ref{contrex2}) une fonction $g$ enti\`ere et transcendante, dont toutes ses d\'eriv\'ees, envoient tout nombre alg\'ebrique $\alpha$ dans $\mathbf{Z}\left[\frac{1}{2},\alpha\right]$. 

\medskip

Afin d'\'etudier l'ensemble $S_{f}$, filtrons l'ensemble d\'enombrable des nombres alg\'ebriques par le degr\'e et la hauteur {\textit {logarithmique}}  absolue (dont la d\'efinition pr\'ecise est rappel\'ee dans la section \ref{notations}): pour $D$ entier $\geq 1$ et $N$ nombre r\'eel $\geq 0$, l'ensemble
$$E_{D,N}= \{ \alpha \in \Qbarre; \espa [\rat (\alpha
):\rat ] \leq D , \esp h(\alpha )\leq N \}$$
est fini et son cardinal $\epsilon_{D,N}$ v\'erifie l'encadrement suivant.

\bigskip 

\begin{lemme}{Encadrement du cardinal de $E_{D,N}.$}\label{e_{D,N}}

Pour tout entier $D \geq 1$ et tout nombre r\'eel $N \geq 0,$ le cardinal $\epsilon _{D,N}$ de $E_{D,N}$ v\'erifie
$$e^{D(D+1)(N-1)} < \epsilon _{D,N} \leq e^{D(D+1)(N+1)}.$$
\end{lemme}
 
 \medskip

Dans ce texte, nous fixons une fonction  $f$ transcendante et nous nous int\'eressons \`a l'ensemble des  $\alpha$ dans $E_{D,N}$ tels que $f(\alpha)$ appartienne \`a $E_{D,N}$.

Le probl\`eme sera local. Pour tout couple de nombres r\'eels $(R,r)$ v\'erifiant $R>r>0$, pour toute fonction $f$ analytique sur $\DRouvert$ et \`a valeurs complexes, pour tout entier $D \geq 1$ et tout nombre r\'eel $N \geq 0,$ on note $\Sigma_{D,N}=\Sigma _{D,N}(f,r)$ l'ensemble des nombres $\alpha$ dans $\Qbarre \cap  \Drferme$ tels que 
$$f(\alpha) \in
\Qbarre ,\esp [\rat(\alpha , f(\alpha )):\rat] \leq D ,\hspace{0,2cm} h(\alpha) \leq N \esp \mathrm{et} \esp h(f(\alpha)) \leq N.$$ 
Ainsi, $S_f$ est la r\'eunion des $\Sigma_{D,N}$ pour $D\geq 1$ et $N\geq 0.$

 
 \medskip

  Dans le cas $D = 1$, J. Pila (th\'eor\`eme 9 de \cite{pila}) et N. Elkies (th\'eor\`eme 4 de \cite{elkies.rat}) obtiennent le r\'esultat suivant. 
 {\textit {Soit $f$  une fonction r\'eelle analytique sur un intervalle ferm\'e, dont l'image n'est pas contenue dans aucune courbe alg\'ebrique. Pour tout $\varepsilon > 0$, il existe une constante $c(f, \varepsilon)$ telle que pour tout $N \geq 1$, le nombre de points rationnels du graphe de $f$ de hauteur inf\'erieure \`a $N$, est inf\'erieur  \`a 
 $$c(f, \varepsilon) \, e^{\varepsilon N}.$$
 }}

 Le r\'esultat suivant montre, en particulier,  que la borne du th\'eor\`eme de J. Pila et de N. Elkies n'est pas loin d'\^etre optimale.

 \begin{thm}\label{contrex}

Soit $\phi$ une fonction positive telle que $\phi(x)/x$ tende vers $0$ quand $x$ tend vers l'infini.  Il existe une suite $(N_{\delta})_{\delta \geq 1}$ de nombres r\'eels, croissant vers l'infini, et une fonction $f$ enti\`ere et transcendante sur $\mathbf{C}(z)$, v\'erifiant
\begin{equation}
\forall \alpha \in \Qbarre, \esp \forall \sigma \geq 0, \esp f^{(\sigma)}(\alpha) \in \rat(\alpha), \label{a)}
\end{equation}
et telle que, pour tout entier $D \geq 1$, pour tout $k\geq D$, 
\begin{equation}
\card(\Sigma _{D,N_{k}}(f,1)) \geq e^{D(D+1)\phi(N_{k})-\log2}. \label{b)}
\end{equation}

\end{thm}


D'autre part, avec les m\'ethodes classiques de transcendance, nous obtenons le r\'esultat suivant.
 
\begin{thm}{Majoration du cardinal de $\Sigma_{D,N}(f,r)$.}\label{propprincipale}

Soient $R$ et $r$ deux nombres r\'eels v\'erifiant $R>r>0$, $ c_0 = \log \left( \frac{R^2 + r^2}{2rR} \right)$ 
et $\delta$ un nombre r\'eel tel que $\delta > 2\left(6/c_{0} \right)^{2}$. 
Soit $f$ une fonction analytique sur $\DRouvert$, continue sur $\DRferme$ et transcendante sur $\mathbf{C}(z)$.

i) Pour tout entier $D \geq1$, il existe une suite de nombres r\'eels $N \geq 0$ tendant vers l'infini pour lesquels
$$\card(\Sigma _{D,N}) < \delta D^3N^2.$$

ii) Pour tout nombre r\'eel $N > 0$, il existe une suite de nombres entiers $D \geq 2 $ tendant vers l'infini pour lesquels
$$\card(\Sigma _{D,N}) < \delta D^3N^2.$$

\end{thm}

%
 Les nombres r\'eels $c_0$ et $\delta$ sont strictement positifs et d\'ependent uniquement de   $r$ et $R$.

 Pour $D$ et $N$ fix\'es, l'ensemble $\Sigma_{D,N}$ est fini, puisqu'il est contenu dans l'ensemble $E_{D,N}$. Le th\'eor\`eme \ref{propprincipale} montre d'une part que, pour $D$ fix\'e et une infinit\'e d'entiers $N$, le cardinal de $\Sigma_{D,N}$ est beaucoup plus petit que $\epsilon_{D,N}$; et, d'autre part que, pour $N$ fix\'e et une infinit\'e d'entiers $D$, la m\^eme conclusion est valable.

Le th\'eor\`eme \ref{contrex} montre qu'on ne peut pas remplacer, dans la partie i) du th\'eor\`eme \ref{propprincipale},  ``il existe une suite de nombres r\'eels $N \geq 0$ tendant vers l'infini" par ``pour tout $N$ assez grand''. De m\^{e}me,  la partie i) du th\'eor\`eme \ref{propprincipale} montre qu'on ne peut pas remplacer, dans le th\'eor\`eme \ref{contrex}, la suite $N_{k}$ par ``pour tout $N$ assez grand''.


\medskip
  
 Le r\'esultat principal de ce travail, le th\'eor\`eme \ref{prop-gen}, g\'en\'eralise le th\'eor\`eme \ref{propprincipale} au cas de plusieurs fonctions analytiques alg\'ebriquement ind\'ependantes. 
  
\medskip

Dans le cas des fonctions r\'eelles, J. Pila (th\'eor\`eme 8 de \cite{pila}) obtient un r\'esultat uniforme \`a la fois en la borne du degr\'e du corps de nombres et en la borne de la hauteur. Cependant, la d\'ependance en le degr\'e n'est pas explicite. On peut d\'eduire de son r\'esultat le corollaire suivant. {\textit {Soit $f$ une fonction analytique sur l'intervalle r\'eel $[-r, r]$. Soient $d$ un entier naturel  non nul et $\varepsilon >0$ un nombre r\'eel.  Il existe un nombre r\'eel $c(f, d, \varepsilon)$ tel que, pour tout corps de nombres $K$ de degr\'e $[K:\rat] = d$ et tout r\'eel $N >0$, le cardinal de l'ensemble des $x$ dans $[-r, r] \cap K$ tels que $f(x) \in K, \esp h(x) \leq N$ et $h(f(x)) \leq N$, est inf\'erieur \`a
$$c(f, d, \varepsilon) \,e^{\varepsilon dN}.$$}}

D'autres r\'esultats sur le cas rationnel, concernant des fonctions r\'eelles, ainsi que des g\'en\'eralisations, ont \'et\'e obtenus pas J. Pila, E. Bombieri et J. Pila, et par J. Pila et A.J. Wilkie. Les m\'ethodes remontent, en partie, \`a celles de \cite{bombieri-pila}. Le plus r\'ecent, \cite{pila-wilkie}, qui porte sur des situations multidimensionnelles, fait intervenir la th\'eorie de la $o$-minimalit\'e.

\bigskip

Cet article est organis\'e de la fa\c con suivante. 
On introduit les notations employ\'ees dans la section \ref{notations}. 
Dans la section \ref{demo-gen}, on \'enonce et on d\'emontre le th\'eor\`eme \ref{prop-gen} qui concerne le nombre de points alg\'ebriques o\`u prennent des valeurs alg\'ebriques plusieurs fonctions analytiques alg\'ebriquement ind\'ependantes. Le th\'eor\`eme \ref{propprincipale} s'en d\'eduit. 
Dans la section \ref{cardinal} on compte le nombre de points alg\'ebriques de hauteur et de degr\'e born\'es; il s'agit du  lemme \ref{e_{D,N}}. On y fait r\'ef\'erence \`a d'autres r\'esultats sur le sujet. 
La construction explicite de la fonction  dont l'existence est assur\'ee dans le th\'eor\`eme \ref{contrex} est donn\'ee dans la section \ref{exemples}; on y donne aussi la construction d'une autre fonction transcendante $g$ envoyant tout nombre alg\'ebrique $\alpha$ dans  $\enteros[\frac{1}{2}, \alpha]$. 

 \bigskip
 

{\small Le lemme \ref{e_{D,N}}, le th\'eor\`eme \ref{contrex}  et la partie i) du th\'eor\`eme \ref{propprincipale} ont fait l'objet d'une annonce dans \cite{macras}. Les d\'etails de leurs preuves sont donn\'es ici. Ces r\'esultats font partie de ma th\`ese de doctorat, r\'ealis\'ee \`a l'Institut de Math\'ematiques de Jussieu. Je remercie Marc Hindry et Michel Waldschmidt de l'avoir dirig\'ee, ainsi que Joseph Oesterl\'e et Fran\c cois Gramain pour leurs commentaires.  }


\section{Notations.}\label{notations}

Pour un corps de nombres $K$, on note $M_K$ l'ensemble de classes d'\'equivalence de ces valeurs absolues  dont la restriction \`a $\mathbf{Q}$ est, soit la valeur absolue archim\'edienne $|.|_{\infty}$, soit une des valeurs absolues ultram\'etriques, normalis\'ees de la fa\c con suivante:

$$|x|_{\infty} = x, \esp \textrm{si} \esp x \in \rat, \esp x > 0,$$ 
$$|p|_p =\frac{1}{p}, \esp \textrm{si} \esp p \esp \textrm{est un nombre premier}.$$ 
La formule du produit pour un \'el\'ement non nul $x$ de $K$ s'\'ecrit alors $\prod_{v \in M_{K}} |x|_v^{d_{v}}=1$, o\`u $d_{v}$ est le degr\'e local en la place $v$.


Si $\alpha $ est un nombre alg\'ebrique et $K$ un corps de nombres le contenant, on d\'efinit sa hauteur projective logarithmique absolue par

$$h(\alpha)= \frac{1}{[K:\rat]} \sum_{v \in M_K} d_{v} \log \max \{ 1, |\alpha|_v \} .$$

Si $\alpha_1, \ldots , \alpha_n$ sont $n$ nombres alg\'ebriques, on a (cf. \cite{w} chapitre 3) les in\'egalit\'es suivantes:
\begin{equation}
h(\alpha_{1} \alpha_{2}) \leq h(\alpha_{1}) + h(\alpha_{2}) \label{h(a.b)<h(a)+h(b)}
\end{equation}
\begin{equation}
h(\alpha_1 + \cdots + \alpha_n ) \leq h(\alpha_1) + \cdots + h(\alpha_n) + \log n.  \label{h(a_1+...+a_d)<h(a_1)+...+h(a_d)}
\end{equation}

Si $a_0X^d+a_1X_{d-1}+\ldots +a_d$ est le polyn\^ome minimal sur
$\mathbf{Z} $ de $\alpha $ (o\`u $d=[\mathbf{Q} (\alpha ):\mathbf{Q} ]$) et $\alpha_1=\alpha, \ldots, \alpha_d$ ses conjugu\'es, on d\'efinit la mesure de Mahler de $\alpha $ par

$$M(\alpha )= \vert a_0\vert \prod_{i=1}^{d}\max \{1, \vert
\alpha_i\vert \}.$$

Elle est reli\'ee \`a la hauteur logarithmique absolue par la
relation (\cite{w}  lemma 3.10)

\begin{equation}
h(\alpha)=\frac{1}{d}\log M(\alpha). \label{mahler_hauteurlog} 
\end{equation}

La hauteur usuelle d'un polyn\^ome $P \in \mathbf{C}\left[ X \right]$ est, par d\'efinition, le maximum du module de ses coefficients. La hauteur usuelle d'un nombre alg\'ebrique $\alpha$ est d\'efinie comme \'etant la hauteur usuelle de son polyn\^ome minimal, \`a savoir (en gardant les m\^emes notations),
$$\mathcal{H}(\alpha)=\max \{|a_0|, \ldots , |a_d|\}.$$
La hauteur usuelle est li\'ee \`a la hauteur logarithmique absolue par la double in\'egalit\'e suivante (\cite{w} chapitre 3 lemma 3.11)

\begin{equation}
\frac {1}{d}\log \mathcal{H}(\alpha )-\log 2 \leq h(\alpha) \leq \frac {1}{d}\log \mathcal{H}(\alpha )+ \frac{1}{2d} \log (d+1). \label{hauteurusuelle_hauteurlog}
\end{equation}
Au lieu de la premi\`ere in\'egalit\'e, nous utiliserons plut\^ot celle-ci:
\begin{equation}
\frac {1}{d}\log \mathcal{H}(\alpha )- \frac{d-1}{d} \log 2 \leq h(\alpha), \label{hauteurusuelle_hauteurloggauche}
\end{equation}
obtenue de la m\^eme fa\c con que la premi\`ere, en remarquant que les coefficients bin\^omiaux sont major\'es par $2^{d-1}.$

Pour un polyn\^ome $P\in \mathbf{C}[X,Y]$ on note $L(P)$ sa longueur. C'est la somme des modules de ses coefficients. 
 
Pour un nombre r\'eel $\rho >0$ et une fonction $F$ continue dans le disque ferm\'e $\Drhoferme,$ on note
$$\vert F\vert _\rho = \max_{\vert z\vert \leq \rho }\vert F(z)\vert.$$
  
Pour $x\in \mathbf{R},$ $[x]$ d\'esigne la partie enti\`ere de $x.$
Elle v\'erifie 
$[x] \in \enteros \hspace{0,2cm}\hbox{et}\hspace{0,2cm}0\leq x-[x] <1$.


\section{Le th\'eor\`eme principal}\label{demo-gen}

Pour des nombres r\'eels $R, \esp r$ et $c_0$ comme dans le th\'eor\`eme \ref{propprincipale} et un entier $t \ge 2,$ on note 
\begin{equation}
\gamma_t = \max_{2 \leq \tau \leq t} \left( \frac{1}{2} \left( \frac{c_0}{3\tau} \right)^{\tau}  \right)^{\frac{-1}{\tau-1}}. \label{gamma}
\end{equation}
Ainsi $\gamma_t$ ne d\'epend que de $R, r$ et $t$.

Pour $f_1, \ldots , f_t$ des fonctions analytiques sur $\DRouvert,$ pour tout entier $D \geq 1$ et tout nombre r\'eel $N \geq 0,$ on consid\`ere l'ensemble $\Sigma _{D,N}(f_1, \ldots , f_t, r)$ des nombres $w \in \complexes \cap \Drferme$ tels que, 
$$\forall i \in \{1, \ldots , t\}, \esp f_i(w) \in \Qbarre, \esp  h(f_i(w)) \leq N $$
et
$$[ \rat (f_1(w), \ldots , f_t(w)): \rat ] \leq D.$$

\begin{thm}{Majoration du cardinal de $\Sigma_{D,N}(f_1, \ldots , f_t,r).$}\label{prop-gen}

Soient $R$ et $r$ des nombres r\'eels v\'erifiant $R>r>0$ et $c_0=\log \left( \frac{R^2 + r^2}{2rR} \right).$  Soient $t$ un entier $\ge 2$ et $\gamma$ une constante r\'eelle telle que $\gamma > \gamma_t.$

Soient $f_1, \ldots , f_t$ des fonctions analytiques sur $\DRouvert$ et continues sur $\DRferme$, alg\'ebri\-quement ind\'e\-pendantes sur $\rat$.

i) Pour tout entier $D \geq1$, il existe une infinit\'e de nombres r\'eels $N \geq 0$ arbitrairement grands pour lesquels
$$\card(\Sigma _{D,N}(f_1, \ldots , f_t,r)) < \gamma\, D^{a(t)}N^{b(t)};$$
\indent ii) pour tout nombre r\'eel $N>0$, il existe une infinit\'e de nombres entiers $D \geq 2$ arbitrairement grands pour lesquels
$$\card(\Sigma _{D,N}(f_1, \ldots , f_t,r)) < \gamma\, D^{a(t)}N^{b(t)},$$
o\`u $a(t) = \frac{t + 1}{t - 1} $ et $b(t) = \frac{t}{t - 1} .$
 
\end{thm}

\bigskip
 
Notons que pour $t \geq 2,$ on a $a(t)\leq 3$ et $b(t) \leq 2.$

\begin{remarque}
Si $t=2$, alors $a(t) = 3$ et $b(t) = 2$. Alors  en prenant pour une des deux fonctions l'identit\'e, l'autre fonction est transcendante, et nous retrouvons le th\'eor\`eme \ref{propprincipale}. 
\end{remarque}

Dans la partie qui suit nous d\'emontrons l'assertion i) du th\'eor\`eme \ref{prop-gen}. L'assertion ii) sera d\'emontr\'e de fa\c con analogue, dans la partie \ref{unif-hauteur}.


\subsection{Le r\'esultat uniforme en la borne du degr\'e.}\label{unif-degre}



La d\'emonstration peut se faire par l'absurde. Voici le sch\'ema de d\'emonstration.

On suppose qu'il existe un entier $D \geq 1$ et un entier $N_0 \geq 1$
suffisamment grand tels que, pour tout $N \geq N_0$, on ait

$$\card(\Sigma _{D,N}(f_1, \ldots, f_t)) \geq \gamma\, D^{a(t)}N^{b(t)}.$$

Pour tout $N\geq N_0,$ on commence par extraire de $\Sigma_{D,N}$ un sous-ensemble convenable $S_{D,N}$ dont on conna\^\i t exactement le nombre d'\'el\'ements.

Le lemme de Siegel nous donne ensuite l'existence d'un polyn\^ome $P$, non
nul, \`a $t$ variables et \`a coefficients entiers (et d\'ependant de
$N_0$), tel que 

$$P(f_1(\omega),\ldots, f_t(\omega))=0,\hspace{0,5cm}  \forall \omega \in S_{D,N_0}.$$

On d\'efinit alors une fonction $F$ sur le disque $\DRferme$ en posant

$$F(z)=P(f_1(z),\ldots,f_t(z)),\hspace{0,5cm}  \forall z \in \DRferme.$$

Ainsi la construction de $P$ entra\^\i ne

$$F(\omega)=0,\hspace{0,5cm} \forall \omega \in S_{D,N_0}.$$

On montrera, gr\^ace \`a une r\'ecurrence, \`a l'in\'egalit\'e de
Liouville et \`a un lemme de Schwarz, qu'on a 

$$F(\omega )=0,\hspace{0,5cm} \forall \omega \in \bigcup _{N \geq N_0}S_{D,N},$$
ce qui impliquera que $F$ est la fonction nulle et contredira le fait que les fonctions $f_{1}, \ldots , f_{t}$ sont alg\'ebriquement ind\'ependantes. D'o\`u le r\'esultat.





Dans cette section, nous noterons $\Sigma_{D,N}$ l'ensemble $\Sigma_{D,N}(f_1, \ldots , f_t,r)$ et $\sigma_{D,N}$ son cardinal.


\subsubsection{Choix des param\`etres et construction de
  $S_{D,N}$.}\label{choix des param}

On fixe un entier $D \ge 1$ et un nombre r\'eel $\gamma > \gamma_t$ et on suppose qu'il existe un nombre r\'eel $N_0$ tel que, pour tout $N \ge N_0,$ 
$$\sigma _{D,N} \geq \gamma \, D^{a(t)}N^{b(t)}.$$

On pose $T = \left[ \frac{c_0 \gamma}{3t} D^{\frac{2}{t-1}} N_0^{\frac{1}{t-1}} \right].$ Comme $b(t)= \frac{t}{t-1}$, quitte \`a augmenter $N_0,$ on a 
$$c_0 \left[ \gamma \, D^{a(t)} (N_0-1)^{b(t)} \right] > $$
\begin{equation}
D\log2 + 2tD\log T + T \left(2t N_0 D + \sum_{i=1}^t \log \max \{1,|f_i|_R \} \right).\label{inegN0}
\end{equation}

On pose $u_1 = \log\left( 2T^{2t} e^{tN_0T} \max \{ 1,|f_1|_R^T \} \ldots \max \{1,|f_t|_R^T \} \right).$ L'in\'egalit\'e (\ref{inegN0}) s'\'ecrit alors:
\begin{equation}
c_0 \left[ \gamma D^{a(t)} (N_0-1)^{b(t)} \right] > u_1 + (D-1) \log(2T^{2t}) + tTN_0(D-1) +
tTN_0D. \label{inegaliteN0}
\end{equation}

Puisque pour tout $N,$ sup\'erieur ou \'egal \`a $N_0,$ le cardinal
$\sigma _{D,N}$ de l'ensemble $\Sigma _{D,N}$ est sup\'erieur ou
\'egal \`a $\gamma D^{a(t)}N^{b(t)},$ on peut extraire de $\Sigma _{D,N}$ un
sous-ensemble $S_{D,N}$ dont le cardinal $s_{D,N}$ est exactement
$\left[ \gamma\, D^{a(t)}N^{b(t)} \right].$ 

Pour $N\geq N_0,$ on num\'erote les \'el\'ements de $S_{D, N} =\{ \omega_1, \omega_2, \ldots ,\omega_{s_{D, N}} \}$.

Comme les nombres $D$ et $N_{0}$ sont fix\'es pour toute la suite, on notera, %
$$S_{D,N_{0}} = S_{0} \esp \textrm{et} \esp s_{D, N_{0}} = s_{0}.$$
%

%


\subsubsection{La fonction auxiliaire.}

\begin{lemme}\label{appli.siegel}

Il existe un polyn\^ome $P \in \enteros [X_1,\ldots, X_t],$ non nul, de
degr\'e en $X_i$ strictement inf\'erieur \`a $T,$ pour tout $i \in \{1, \ldots, t\},$
tel que 
$$P(f_1(\omega),\ldots, f_t(\omega))=0, \hspace{0,5cm} \forall \omega \in S_{0}$$
et dont les coefficients sont major\'es en valeur absolue par $2T^t
e^{tTN_0}.$

\end{lemme}

\noindent \textbf{D\'emonstration du lemme \ref{appli.siegel}.}

On \'ecrit le polyn\^ome cherch\'e sous la forme
$$P(X_1,\ldots ,X_t)=\sum_{i_1=0}^{T-1} \ldots \sum_{i_t=0}^{T-1}c_{i_1 \ldots i_t} \, X_1^{i_1} \ldots X_t^{i_t}.$$
Montrons que 
$$\sum_{i_1=0}^{T-1} \ldots \sum_{i_t=0}^{T-1}c_{i_1, \ldots, i_t} \, f_1(\omega_k)^{i_1} \ldots f_t(\omega_k)^{i_t} = 0,
\hspace{0,5cm} \forall k\in \{1,\ldots , s_{0}\}$$ 
poss\'ede une solution non triviale $c_{i_1, \ldots, i_t}$ dans $\enteros^{T^{t}}$. Pour ainsi faire, on applique le lemme 1.1 de \cite{gmw} en posant : 
$$\tau = t, \esp N_1 = \ldots = N_t = T,  \esp \mu=s_{0}, \esp L=T^t,$$
$$\{ \alpha_{h,1}, \ldots, \alpha_{h,{\tau}} \}_{h=1, \ldots, \mu} = \{f_1(\omega_k), \ldots, f_t(\omega_k)\}_{k=1, \ldots ,
  s_{0}},$$
et, pour $k \in \{1,\ldots, s_{0} \}$, 
$d_k=[\rat(f_1(\omega_k), \ldots, f_t(\omega_k)):\rat]$.
Comme $\gamma > \gamma_t,$ et $t \ge 2,$ d'apr\`es (\ref{gamma}), nous avons $2\gamma < \left(\frac{c_0 \gamma}{3t} \right)^t $ et, quitte \`a augmenter $N_{0}$,  $T$ v\'erifie 
\begin{equation}
2\, \gamma \, D^{\frac{2t}{t-1}}N_0^{\frac{t}{t-1}} < T^t \le \left( \frac{c_0 \gamma}{3t} \right)^t D^{\frac{2t}{t-1}} \, N_0^{\frac{t}{t-1}}. \label{double T}
\end{equation}
Alors
$$\sum_{h=1}^{\mu}d_h = \sum_{k=1}^{s_{0}}d_k \leq s_{0}D \leq \gamma \, D^{\frac{2t}{t-1}} \, N_0^{\frac{t}{t-1}} < T^t,$$
et l'hypoth\`ese $L>\sum_{h=1}^{\mu}d_h$ est satisfaite.

Le syst\`eme poss\`ede donc une solution
$c=(c_{i_1, \ldots, i_t})_{0\leq i_1, \ldots, i_t\leq T-1} \, ,$
o\` u les $c_{i_1, \ldots, i_t}$ sont des entiers rationnels non tous nuls qui
v\'erifient
$$\max _{0\leq i_1, \ldots, i_t\leq T-1}\vert c_{i_1, \ldots, i_t} \vert \leq \left[ \left( 2^{s_{0}} \prod_{k=1}^{s_{0}}M_k \right)^{\frac{1}{T^t-s_{0}D}} \right]$$
o\`u $M_k = T^{td_k} \prod_{r=1}^t H(f_r(\omega_k))^{(T-1)d_k}$.
Comme pour tout $k \in  \{1, \ldots, s_{0}\}$  on a $d_k \leq D$ et $\omega_k \in S_{D, N_{0}}$, alors pour tout $r \in \{1, \ldots, t\}$, on a $H(f_r(\omega_k)) \leq e^{N_0}$. D'o\`u
$$M_k\leq T^{tD}e^{tN_0(T-1)D}$$
et
$\max _{0\leq i_1, \ldots, i_t\leq T-1}\vert c_{i_1, \ldots, i_t} \vert \leq \left(
  2^{s_{0}}T^{t D s_{0}}e^{t N_0 T D s_{0}}
\right)^{\frac{1}{T^t-s_{0}D}}$.
Or $D\geq 1,\esp s_{0}\geq 1$ et, d'apr\`es (\ref{double T}), on a $T^t > 2s_{0}D$, ce qui
conduit \`a la majoration 
$$\max _{0\leq i_1,i_t\leq T-1}\vert c_{i_1, \ldots, i_t} \vert \leq 2\, T^t \, e^{tTN_0}$$
et d\'emontre le lemme. \hfill $\Box$

\bigskip 

On d\'efinit une fonction $F$ continue sur le disque ferm\'e $\DRferme$,  en posant
$$F(z)=P(f_1(z), \ldots, f_t(z)),\hspace{0,5cm}\forall z \in \DRferme.$$
%


\subsubsection{La r\'ecurrence.}

On va montrer, par r\'ecurrence, que 
$$\forall N \geq N_0, \esp \forall \omega \in S_{D,N}, \esp F(\omega )=0.$$
Par construction du polyn\^ome $P$, on a $F(\omega)=0$ pour tout $\omega$ appartenant \`a $S_{D, N_0}$, et les deux lemmes suivants permettront de conclure.

\medskip

\begin{lemme}\label{A}
Pour tout $N\geq N_0,$ on a
$$\left[ \forall  \omega \in S_{D,N}, \esp F(\omega)=0 \right] \Rightarrow \left[ \forall  \omega \in S_{D,N+1}, \esp  |F(\omega)| \leq e^{u_1-c_0 s_{D, N}} \right].$$
\end{lemme}

\begin{lemme}\label{B}
Pour tout $N\geq N_0 + 1$ et tout $\omega \in S_{D,N},$ on a
$$|F(\omega)| \leq e^{u_1-c_0 s_{D, N-1}} \Rightarrow F(\omega)=0.$$
\end{lemme}

\bigskip

\noindent \textbf{D\'emonstration du lemme \ref{A}.}

On se donne $N\geq N_0$ et on suppose que $F(\omega)=0$ pour tout $\omega \in S_{D, N}$. Nous appliquons un lemme de Schwarz (cf. \cite{w} exercice 4.3) \`a la fonction $F$,  continue sur $\DRferme$, holomorphe sur $\DRouvert$, s'annulant sur chaque point $\zeta_i=\omega_i$ du disque  $\DRferme$ avec une multiplicit\'e $\geq \sigma_{i} = 1$.
On obtient 
$$\vert F \vert_r \leq \vert F\vert _R\prod_{k=1}^{s_{D, N}}\left(
  \frac{R^2+r\vert \omega_k\vert}{R(r+\vert\omega_k\vert)}\right)^{-1}.$$
Or $\frac{R^2+r\vert
  \omega_k\vert}{R(r+\vert\omega_k\vert)}-\frac{R^2+r^2}{2rR}
=\frac{(R^2-r^2)(r-\vert\omega_k\vert)}{2rR(r+\vert
  \omega_k\vert)}\geq 0 $ d\`es que $\vert
  \omega_k\vert \leq r$, 
donc, pour $\omega_k \in S_{D, N}$ on a
$$\frac{R^2+r\vert \omega_k\vert}{R(r+\vert\omega_k\vert)} \geq
\frac{R^2+r^2}{2rR},$$
et comme $c_0 = \log \left( \frac{R^2 + r^2}{2rR} \right)$, on a
$\vert F \vert_r \leq \vert F\vert _R \,e^{-c_0s_{D, N}}$. 
De plus, pour tout $z$ dans $\DRferme$, on a
$$F(z)=P(f_1(z), \ldots,f_t(z))=\sum_{i_1=0}^{T-1} \ldots \sum_{i_t=0}^{T-1}c_{i_1, \ldots, i_t} \,f_1(z)^{i_1} \ldots f_t(z)^{i_t}.$$
On peut donc majorer le module de $F(z)$ en s'aidant de la borne
obtenue pour les coefficients $c_{i_1, \ldots, i_t}$ au lemme
\ref{appli.siegel}:
$$\vert F(z)\vert   \leq   T^t\max_{0\leq i_1, \ldots, i_t\leq T-1} \{|c_{i_1, \ldots, i_t}|\} \,\max \{ 1,|f_1|_R^T \} \ldots \max \{1,|f_t|_R^T \}$$
$$ \hspace{2,2cm}      \leq   2 \,T^{2t} \,e^{tTN_0} \max \{ 1,|f_1|_R^T \} \ldots \max \{1,|f_t|_R^T \}.$$
Comme $u_1 = \log \left( 2\,T^{2t} \,e^{tTN_0} |f_1|_R^T\ldots |f_t|_R^T \right),$ alors $\vert F\vert _R\leq  e^{u_1}$ et
$$\max_{|z|\leq r}\vert F(z)\vert=\vert F\vert_r\leq e^{u_1-c_0s_{D, N}}.$$
En particulier, tout $\omega \in S_{D, N+1}$ est de module $\leq r,$ donc $\vert F(\omega)\vert \leq e^{u_1-c_0s_{D, N}},$ ce qui d\'emontre le lemme \ref{A}.
\hfill $\Box$

\bigskip

\bigskip

\noindent \textbf{D\'emonstration du lemme \ref{B}.}

On se donne $N\geq N_0 + 1, \esp \omega $ dans $S_{D, N}$, et on suppose que  $\vert F(\omega)\vert \leq e^{u_1-c_0s_{N-1}}$.

Comme $s_{D, N-1} = [\gamma D^{a(t)} (N-1)^{b(t)}]$, et qu'au paragraphe \ref{choix des param} nous avons suppos\'e que $N_0$ v\'erifiait (\ref{inegaliteN0}), alors $N$ v\'erifie:
\begin{equation}
c_0 s_{D, N-1} > u_1 + (D-1) \log(2T^{2t}) + tTN_0(D-1) +
tTND. \label{inegaliteN}
\end{equation}

D'autre part, le lemme \ref{appli.siegel} nous a donn\'e un majorant pour la
longueur du polyn\^ome $P$, \`a savoir,
$$L(P)=\sum_{i_1=1}^{T-1} \ldots \sum_{i_t=1}^{T-1} |c_{i_1, \ldots, i_t}| \leq T^t \max_{0\leq i_1, \ldots, i_t\leq T-1} \{|c_{i_1, \ldots, i_t}|\} \leq 2T^{2t}e^{tTN_0},$$
et comme $\omega$ appartient \`a $S_{D, N}$, on a $H(f_r(\omega)) \leq e^N$, pour tout $r$ dans $ \{1, \ldots, t\}$. Donc l'in\'egalit\'e (\ref{inegaliteN}) donne
$$e^{u_1-c_0s_{D, N-1}}< L(P)^{-(D-1)} \prod_{r=1}^t H(f_r(\omega))^{-DT}.$$

D'apr\`es l'hypoth\`ese, on en d\'eduit
$$\vert F(\omega)\vert < L(P)^{-(D-1)} \prod_{r=1}^t H(f_r(\omega
))^{-DT}.$$

Pour conclure que $F(\omega)=0,$ nous appliquons l'in\'egalit\'e de Liouville (proposition 3.14 de \cite{w}), ou plut\^ot sa contrapos\'ee,  en posant: 
$$K=\mathbf{Q}(f_1(\omega), \ldots, f_t(\omega)), \hspace{0,2cm}\textrm{corps
  de nombres de degr\'e}\leq D,\hspace{0,2cm} |.|_v=|.|,$$
$$l=t, \esp N_1= \ldots = N_t = T, \esp \gamma=(f_1(\omega), \ldots, f_t(\omega)) \esp \textrm{et}\esp f=P.$$
%
%
%

On en conclut que $P(f_1(\omega), \ldots,  f_t(\omega))=0, \hspace{0,2cm}\mathrm{i.e.}
\hspace{0,2cm} F(\omega)=0,$ ce qui d\'emontre le lemme \ref{B} et donc la r\'ecurrence. \hfill $\Box$


\subsubsection{Conclusion.}

Nous avons montr\'e que pour tout $N\geq N_0,$ la fonction $F$ s'annule sur
l'ensemble $S_{D, N}$ et donc
\begin{equation}
F(\omega )=0, \hspace{0,5cm} \forall \omega \in \bigcup_{N\geq N_0}
S_{D, N}. \label{concl}
\end{equation}

Les ensembles $S_{D, N}$ sont tous inclus dans le compact
$\Drferme$ et la fonction $F$ est holomorphe sur l'ouvert $\DRouvert$
qui contient $\Drferme.$
De plus, $F$ n'est pas la fonction identiquement nulle car elle est
d\'efinie comme \'etant la valeur de $P \in \enteros[X_1,\ldots,X_t]$ au point
$(f_1(z),\ldots,f_t(z))$ o\`u les $f_i$ sont des fonctions alg\'ebriquement ind\'ependantes sur $\bf{Q}.$
Ainsi $F$ ne peut pas s'annuler sur une infinit\'e de points de
$\Drferme,$ et donc en particulier sur $\bigcup_{N\geq N_0}S_{D, N}.$

Ceci montre que  (\ref{concl}) contredit le th\'eor\`eme sur les z\'eros isol\'es d'une
fonction holomorphe et nous donne le r\'esultat cherch\'e. \hfill $\Box$


\subsection{Le r\'esultat uniforme en la borne de la hauteur.}\label{unif-hauteur}

L'assertion ii) du th\'eor\`eme \ref{prop-gen} se d\'emontre de fa\c con analogue \`a l'assertion i). 


\subsubsection{Choix des param\`etres et construction de
  $S_{D', N_0}$.}\label{choix des param-sym}

On fixe un nombre r\'eel $N_0$ et un nombre r\'eel $\gamma > \gamma_t$ et on suppose qu'il existe un nombre entier $D$ tel que, pour tout $D'\geq D$, on ait $\sigma _{D', N_0} \geq \gamma \, D'^{a(t)}N_{0}^{b(t)}$.
On pose, comme auparavant,  $T = \left[ \frac{c_0 \gamma}{3t} D^{\frac{2}{t-1}} N_0^{\frac{1}{t-1}} \right].$ Comme $a(t)= \frac{t+1}{t-1}$, quitte \`a augmenter $D$, on a 
$$c_0 \left[ \gamma \, (D-1)^{a(t)} N_0^{b(t)} \right] > $$
\begin{equation}
(D+1)\,\log(2T^{2t}) + 2t N_{0} T (D+1) + T \sum_{i=1}^t \log \max \{1,|f_i|_R \}.\label{inegD}
\end{equation}
On pose encore  $u_1 = \log\left( 2T^{2t} e^{tN_0T} \max \{ 1,|f_1|_R^T \} \ldots \max \{1,|f_t|_R^T \} \right)$. L'in\'egalit\'e (\ref{inegD}) s'\'ecrit alors:
\begin{equation}
c_0 \left[ \gamma (D-1)^{a(t)} N_0^{b(t)} \right] > u_1 + D \, \log(2T^{2t}) + tN_0T(2D +1), \label{inegaliteD}
\end{equation}
et vient remplacer l'in\'egalit\'e (\ref{inegaliteN}) de la section \ref{unif-degre}.

Puisque pour tout $D'$, sup\'erieur ou \'egal \`a $D$, le cardinal
$\sigma _{D', N_0}$ de l'ensemble $\Sigma _{D', N_0}$ est sup\'erieur ou
\'egal \`a $\gamma D'^{a(t)}N_{0}^{b(t)},$ on peut extraire de $\Sigma _{D', N_0}$ un
sous-ensemble $S_{D', N_0}$ dont le cardinal $s_{D', N_0}$ est exactement
$\left[ \gamma\, D'^{a(t)}N_{0}^{b(t)} \right].$ 

Comme les nombres $N_{0}$ et $D$ sont fix\'es pour toute la suite, on notera, 
$$ S_{D, N_0}=S_0 \esp \textrm{et} \esp s_{D, N_0}=s_0.$$
Pour $D'\geq D$, on notera : $S_{D', N_{0}}=\{ \omega_1, \omega_2, \ldots ,\omega_{s_{D', N_{0}}} \}$.
%


\subsubsection{La fonction auxiliaire.}

Le lemme \ref{appli.siegel} s'applique encore ici et on d\'efinit une fonction $F$ continue sur le disque ferm\'e $\DRferme$,  en posant
$$F(z)=P(f_1(z), \ldots, f_t(z)),\hspace{0,5cm}\forall z \in \DRferme.$$
%


\subsubsection{La r\'ecurrence.}

On va montrer, par r\'ecurrence, que 
$$\forall D' \geq D, \esp \forall \omega \in S_{D', N_0}, \esp F(\omega )=0.$$
Par construction du polyn\^ome $P$, on a $F(\omega)=0$ pour tout $\omega$ dans $S_{D, N_{0}},$ et les deux lemmes suivants permettent de conclure.

\medskip

\begin{lemme}\label{Asym}
Pour tout $D'\geq D$, on a
$$\left[ \forall  \omega \in S_{D', N_0}, \esp F(\omega)=0 \right] \Rightarrow \left[ \forall  \omega \in S_{D' +1, N_0}, \esp  |F(\omega)| \leq e^{u_1-c_0 s_{D', N_{0}}} \right].$$
\end{lemme}

\begin{lemme}\label{Bsym}
Pour tout $D'\geq D$ et tout $\omega \in S_{D', N_0},$ on a
$$|F(\omega)| \leq e^{u_1-c_0 s_{D'-1, N_{0}}} \Rightarrow F(\omega)=0.$$
\end{lemme}

\bigskip

\noindent \textbf{D\'emonstration du lemme \ref{Asym}.}

Avec les notations choisies, la d\'emonstration est la m\^eme que celle du lemme \ref{A}.

\bigskip
\bigskip

\noindent \textbf{D\'emonstration du lemme \ref{Bsym}.}

On se donne $D' \geq D + 1, \esp \omega$ dans $S_{D', N_0},$ et on suppose que  $\vert F(\omega)\vert \leq e^{u_1-c_0 \,s_{D'-1, N_{0}}}$. Par rapport \`a la section \ref{unif-degre}, il faut remarquer qu'ici on a $[\rat(\omega):\rat] \leq D'$ et $h(\omega) \leq N_{0}$.

Comme $s_{D'-1, N_{0}} = [\gamma (D'-1)^{a(t)} N_{0}^{b(t)}],$ et qu'au paragraphe \ref{choix des param-sym} nous avons suppos\'e que $D$ v\'erifiait (\ref{inegaliteD}), alors $D'$ v\'erifie:
\begin{equation}
c_0 s_{D'-1, N_{0}} > u_1 + D \, \log(2T^{2t}) + tTN_0(2D +1). \label{inegaliteN}
\end{equation}

D'autre part, le lemme \ref{appli.siegel} nous a donn\'e un majorant pour la
longueur du polyn\^ome $P$, \`a savoir, $L(P) \leq 2T^{2t}e^{tTN_0}$, 
et comme $\omega$ appartient \`a $S_{D', N_0}$, on a $H(f_r(\omega)) \leq e^{N_{0}}$, pour tout $r$ dans $ \{1, \ldots, t\}$, et donc l'in\'egalit\'e (\ref{inegaliteN}) donne
$$e^{u_1-c_0 \,s_{D'-1, N_{0}}}< L(P)^{-D} \prod_{r=1}^t H(f_r(\omega))^{-(D+1)T}.$$

D'apr\`es l'hypoth\`ese, on a donc
$$\vert F(\omega)\vert < L(P)^{-D} \prod_{r=1}^t H(f_r(\omega
))^{-(D+1)T}.$$

Pour conclure que $F(\omega)=0$, nous appliquons la contrapos\'ee de l'in\'egalit\'e de Liouville (proposition 3.14 de \cite{w}),  en posant: 
$$K=\mathbf{Q}(f_1(\omega), \ldots, f_t(\omega)), \hspace{0,2cm}\textrm{corps
  de nombres de degr\'e}\leq D +1,\hspace{0,2cm} |.|_v=|.|,$$
$$l=t, \esp N_1= \ldots = N_t = T, \esp \gamma=(f_1(\omega), \ldots, f_t(\omega)) \esp \textrm{et}\esp f=P.$$
%
%
%

On en conclut que $P(f_1(\omega), \ldots,  f_t(\omega))=0, \hspace{0,2cm}\mathrm{i.e.}
\hspace{0,2cm} F(\omega)=0,$ ce qui d\'emontre le lemme \ref{B} et donc la r\'ecurrence. \hfill $\Box$


\subsubsection{Conclusion.}

On conclut de la m\^eme fa\c con que dans la section \ref{unif-degre}, \`a l'aide du th\'eor\`eme sur les z\'eros isol\'es d'une fonction holomorphe. \hfill $\Box$


\section {Estimation du cardinal de $E_{D,N}$}\label{cardinal}

On rappelle que pour $D$ entier $\geq 1$ et $N$ r\'eel $\geq 0,$ on a not\'e $E_{D,N}$ l'ensemble des nombres alg\'ebriques de degr\'e major\'e par $D$ et de hauteur logarithmique absolue major\'ee par $N.$

Pour estimer son cardinal $\epsilon_{D,N},$ nous consid\'erons l'ensemble $A_{D,H}$ ($H$ \'etant un nombre r\'eel positif) des nombres alg\'ebriques de degr\'e major\'e par $D$ et de hauteur usuelle major\'ee par $H,$
$$A_{D,H}= \{ \alpha \in \Qbarre ; \hspace{0,2cm}[\rat (\alpha
):\rat ] \leq D , \hspace{0,2cm}\mathcal{H}(\alpha )\leq H \}.$$

Pour tout entier naturel $d \geq 1$ nous consid\'erons aussi l'ensemble $\mathcal{A}_{d,H}$ des nombres alg\'ebriques de degr\'e exactement $d$ et de hauteur usuelle major\'ee par $H$,  
$$\mathcal{A}_{d,H}= \{ \alpha \in \Qbarre ; \hspace{0,2cm}[\rat (\alpha
):\rat ] =d , \hspace{0,2cm}\mathcal{H}(\alpha )\leq H \},$$
et l'ensemble $\mathcal{P} _{d,H}$ des polyn\^omes \`a coefficients entiers, \`a une variable, non nuls, de degr\'e exactement $d$, de hauteur usuelle major\'ee par $H,$ irr\'eductibles dans $\mathbf{Z}[X]$ (donc les coefficients sont premiers entre eux dans leur ensemble) et dont le coefficient dominant est positif. Un polyn\^ome $P$ de $\mathcal{P} _{d,H}$ s'\'ecrit:
$$P(X)=a_0X^d+ a_1X^{d-1}+\cdots + a_d$$ 
avec $a_i \in \mathbf{Z},\hspace{0,2cm} |a_i|\leq H, \hspace{0,2cm} \forall i\in\{1, \ldots, d\},\hspace{0,2cm} \pgcd(a_0, \ldots , a_d)=1, \hspace{0,2cm} 1\leq a_0 \leq H.$

Pour tout \'el\'ement $\alpha$ de $A_{D,H}$, il existe un entier $d\leq D$, tel que $\alpha $ appartienne \`a $\mathcal{A}_{d,H}$, et tout $\alpha$ de $\mathcal{A}_{d,H}$ est un z\'ero d'un polyn\^ome $P$ de $\mathcal{P}_{d,H}$. (Les $d-1$ autres racines -distinctes-  de $P$ sont les conjugu\'es de $\alpha$.) Inversement, tout z\'ero d'un polyn\^ome de $\mathcal{P}_{d,H}$ est un \'el\'ement de $\mathcal{A}_{d,H}.$ On remarque ainsi que l'ensemble $A_{D,H}$ est l'union disjointe des ensembles $\mathcal{A}_{d,H}$, pour ${1\leq d\leq D}$. Donc
\begin{equation}
\card(A_{D,H})=\sum_{d=1}^D \card(\mathcal{A}_{d,H}) \label{aster}
\end{equation}
et
\begin{equation}
\card(\mathcal{A}_{d,H}) = d \, \card(\mathcal{P} _{d,H}) \, , \hspace{0,3cm} \forall d \geq 1. \label{aster2}
\end{equation}


\subsection{Estimation du cardinal de $A_{D, H}$.}

\begin{lemme}{Encadrement du cardinal de $A_{D,H}$.}\label{card A_{D,H}}

Pour tout entier $D \geq 1$ et tout nombre r\'eel $H \geq 1,$ on a 
$$ \card(A_{D,H}) \leq DH(2H+1)^D;$$
pour tout entier $D \geq 2$ et tout nombre r\'eel $H \geq 2,$
$$\frac{D}{8}(H-2)^{D+1} < \card(\mathcal{A}_{D,H}) \leq \card(A_{D,H})$$
et pour $D=1$ et tout nombre r\'eel $H\geq 1,$
$$(H-1)^2 < \card(\mathcal{A}_{1,H}) = \card(A_{1,H}). $$

\end{lemme}

Nous allons maintenant d\'emontrer le lemme \ref{card A_{D,H}}.


\subsubsection{Majoration de $\card(A_{D,H})$.}

On fixe $D\geq 1$ et $H \geq 1.$
Comme l'union $\bigcup_{1\leq d\leq D}\mathcal{P}_{d,H}$ est disjointe, en utilisant (\ref{aster}) et (\ref{aster2}), on a
$$\card(A_{D,H}) = \sum_{d=1}^D d \, \card(\mathcal{P}_{d,H}) \leq D\sum_{d=1}^D \card(\mathcal{P}_{d,H}) = D \, \card ( \bigcup_{1\leq d\leq D} \mathcal{P}_{d,H}).$$
Or le cardinal de $\bigcup_{1\leq d\leq D} \mathcal{P}_{d,H}$ est inf\'erieur au nombre de polyn\^omes non nuls $Q\in \mathbf{Z}[X]$ qui s'\'ecrivent
$$Q(X)=a_0X^D+ a_1X^{D-1}+\cdots + a_D$$
avec $|a_i|\leq H$ pour tout $i\in\{1, \ldots, D\},$ de coefficient dominant $a_{0}$ strictement positif et inf\'erieur ou \'egal \`a $H$. En effet, au polyn\^ome $P(X) = a_0 X^d + \cdots + a_d$ de $\bigcup_{1\leq d\leq D} \mathcal{P}_{d,H}$ on peut associer le polyn\^ome $Q(X) = P(X) X^{D-d }$.

Comme il existe $[H]$ nombres entiers $>0$  de valeur absolue inf\'erieure ou \'egale \`a $H,$ il y a autant de choix pour le coefficient dominant, et $2[H]+1$ pour les $D$ autres coefficients (\'eventuellement nuls). On en d\'eduit qu'il existe $[H](2[H]+1)^D$  tels polyn\^omes $Q$, et
$$\card(A_{D,H}) \leq DH(2H+1)^D.$$


\subsubsection{Minoration de $\card(A_{D,H})$.}

Soient $D$ un entier $\geq 1$ et $H$ un nombre r\'eel $\geq 0.$ 
D'apr\`es (\ref{aster}),  $\card (A_{D,H}) \geq \card(\mathcal{A}_{D,H}) $,  et d'apr\`es (\ref{aster2}), 
%
%
 pour minorer le cardinal de $A_{D,H},$ il suffit de minorer celui de $\mathcal{P}_{D,H}$.

Pour $D \geq 2$ et $H \geq 2,$ au lieu de consid\'erer l'ensemble $\mathcal{P}_{D,H},$ nous pouvons juste compter les polyn\^omes $P$ de $\mathcal{P}_{D,H}$ tels que $\pgcd (a_0,a_1)=1$. La condition $\pgcd(a_0, \ldots , a_D)=1$ est alors automatiquement v\'erifi\'ee. Parmi ces polyn\^omes-l\`a, nous pouvons encore nous restreindre \`a ceux qui v\'erifient le crit\`ere d'irr\'eductibilit\'e d'Eisenstein (cf. \cite{lang.algebra}) sur l'anneau des entiers $\enteros$, pour le nombre premier 2.  
On rappelle que le polyn\^ome $P$ est 2-Eisenstein sur $\enteros$ (et donc irr\'eductible sur $\rat$), si
$$2 \hspace{-0,2cm}\not| a_0, \hspace{0,2cm}  2| a_i, \hspace{0,2cm} \forall i\in \{1, \ldots , D\}\hspace{0,2cm} \textrm{et} \hspace{0,2cm} 4 \hspace{-0,2cm}\not| a_D.$$

On note $\mathcal{E}_{D,H}(2)$ l'ensemble de ces polyn\^omes $P$. Comme il est inclus dans $ \mathcal{P}_{D,H}$, on a 
$$\card(\mathcal{A}_{D,H})  \geq D \, \card(\mathcal{E}_{D,H}(2))$$
et nous sommes ramen\'es \`a minorer le cardinal de  $\mathcal{E}_{D,H}(2)$.

On a
$$\card\{a\in \mathbf{Z}/ \hspace{0,2cm} |a|\leq H, \hspace{0,2cm} 2|a\}= 2\left[\frac{H}{2}\right]+1 > H-1,$$
et
$$\card\{a\in \mathbf{Z}/ \hspace{0,2cm} |a|\leq H, \hspace{0,2cm} 2|a,\hspace{0,2cm} 4 \hspace{-0,2cm}\not| a \}= 2\left[ \frac{H+2}{4}\right] > \frac{H-2}{2}.$$
En notant $c_H$ le nombre de couples $(a_0,a_1) \in \enteros ^2 $ tels que 
$$1\leq a_0\leq H, \hspace{0,2cm} |a_1| \leq H, \hspace{0,2cm} 2|a_1, \hspace{0,2cm} \pgcd(a_0, a_1)=1,$$
nous obtenons
$$\card(\mathcal{E}_{D,H}(2)) = c_H \left(2\left[\frac{H}{2}\right]+1\right)^{D-2}\left(2\left[ \frac{H+2}{4}\right]\right),$$
et donc (comme $H \geq 2$), une minoration de ce cardinal, en fonction de $c_H,$
$$\card(\mathcal{E}_{D,H}(2)) > c_H (H-2)^{D-2} \left(\frac{H-2}{2}\right) = \frac{c_H}{2} (H-2)^{D-1}.$$
Il existe $[\frac{H+1}{2}]$ nombres impairs compris entre 1 et $H$, et $[\frac{H}{2}]$ nombres pairs compris entre 1 et $H.$ En tenant compte du signe de $a_1,$ cela implique que le cardinal de l'ensemble 
$$\{(a_0,a_1)\in \enteros ^2 / \hspace{0,2cm} 1\leq a_0\leq H, \hspace{0,2cm} 0<|a_1| \leq H, \hspace{0,2cm} 2 \hspace{-0,2cm}\not|a_0, \hspace{0,2cm} 2|a_1\}$$
est \'egal \`a $2\left[\frac{H+1}{2}\right] \left[\frac{H}{2}\right]$,
et nous avons $2\left[\frac{H+1}{2}\right] \left[\frac{H}{2}\right] \geq \frac{H(H-2)}{2}$ (voir les deux cas, pour $n$ entier positif, $2n \leq H < 2n+1$ et $2n+1 \leq H < 2n+2).$

Pour minorer $c_H,$ on tient compte du couple $(1,0)$ et on enl\`eve tous les couples $(a_0, a_1)$ dont le $\pgcd$ est divisible par un nombre premier $p$ (on a donc $p\geq 3$ car $2 \hspace{-0,2cm} \not| a_0.)$ Ici, $p$ d\'esignera toujours un nombre premier. Or, pour un nombre premier $p\geq 3$, on a
$$\card \{a\in \enteros  / \hspace{0,2cm} 1\leq a\leq H, \hspace{0,2cm} 2 \hspace{-0,2cm}\not|a, \hspace{0,2cm} p|a\} \leq \card \{a\in \enteros  / \hspace{0,2cm} 1\leq a\leq H, \hspace{0,2cm} p|a\}  = \left[\frac{H}{p}\right]$$
et $\card \{a\in \enteros \setminus \{0\}  / \hspace{0,2cm} |a|\leq H, \hspace{0,2cm} 2|a, \hspace{0,2cm} p|a\}$ est \'egal au nombre de $a\in \enteros \setminus \{0\}$ tels que $|a|\leq H$ et $2p|a$, c'est-\`a-dire $2\left[\frac{H}{2p}\right]$ d'o\`u,
$$c_H  \geq 1+ \frac{H(H-2)}{2}  -\sum_{p\geq 3} 2\left[\frac{H}{p}\right]\left[\frac{H}{2p}\right]  \geq 1+\frac{H(H-2)}{2} -H^2 \sum_{p\geq 3}\frac{1}{p^2}.$$
Comme, pour tout nombre $x\in [0,1[,$ on a $-\log (1-x)= \sum_{ k\geq 1} \frac{x^k}{k} \geq x,$ et que, d'une part, $\zeta(2) = \sum_{k \geq 1} \frac{1}{k^2} = \prod_{p \geq 2} \left(\frac{1}{1-\frac{1}{p^2}}\right),$ et d'autre part, $\zeta(2) = \frac{\pi ^2}{6},$ alors
$$\sum_{p\geq 2} \frac{1}{p^2} \leq \sum_{p\geq 2} -\log \left(1-\frac{1}{p^2}\right) = \log (\zeta (2)) = \log \left(\frac{\pi^2}{6} \right) \leq \frac{1}{2}, \hspace{0,2cm} \textrm{d'o\`u} \hspace{0,2cm} \sum_{p\geq 3 } \frac{1}{p^2} \leq \frac{1}{4},$$
et ainsi
$$c_H \geq 1+\frac{H(H-2)}{2} -\frac{H^2}{4} = \frac{(H-2)^2}{4}.$$
On en d\'eduit 
$$\card(\mathcal{E}_{D,H}(2)) > \frac{1}{8}(H-2)^{D+1},$$
et finalement, pour tout $D \geq 2$ et tout $H\geq 2$,  
$\card(\mathcal{A}_{D,H}) > \frac{D}{8}(H-2)^{D+1}$.

Pour $D=1$ et tout nombre r\'eel $H \geq 1,$ on a $ A_{1,H} = \mathcal{A}_{1,H}$ et leur cardinal est \'egal \`a celui de l'ensemble des polyn\^omes $P\in \mathbf{Z}[ X ]$ qui s'\'ecrivent
$$P(X)=aX+b \hspace{0,2cm} \textrm{avec} \hspace{0,2cm} 1\leq a \leq H, \hspace{0,2cm} |b| \leq H , \hspace{0,2cm} \pgcd(a,b)=1.$$
Avec des calculs analogues \`a ceux faits pr\'ec\'edemment, on trouve
$$\card (A_{1,H}) \geq  2\left[ H \right]^2 - \sum_{p\geq 2}\left( \left[\frac{H}{p}\right] 2 \left[\frac{H}{p}\right] \right)$$
et alors
$$\card (A_{1,H}) \geq 2\left[H\right]^2 - 2\left[H\right]^2 \sum_{p \geq 2} \frac{1}{p^2} \geq 2\left[H\right]^2 - 2\left[H\right]^2 \frac{1}{2} = \left[H\right]^2  > (H-1)^2. $$
\hfill $\Box$


\subsection{D\'emonstration du lemme \ref{e_{D,N}}.}

%
%

\subsubsection{Majoration de $\card(E_{D,N})$.}

Soient $D$ un entier $\geq 1$ et $N$ un nombre r\'eel positif ou nul.
Soit $\alpha$ dans $E_{D,N}$ avec $\deg(\alpha )=d.$ Comme $d\leq D$ et $h(\alpha)\leq N$, d'apr\`es (\ref{hauteurusuelle_hauteurloggauche}),
%
%
$$\mathcal{H}(\alpha )\leq 2^{d-1} \, e^{d\, h(\alpha)} \leq 2^{D-1}\, e^{DN}.$$

On a donc $E_{D,N} \subset A_{D,H}$ pour $H=2^{D-1}\, e^{DN}$  et, en utilisant la majoration obtenue pour le cardinal de $A_{D,H}$ au lemme \ref{card A_{D,H}} (on a $H \geq 1$ car $D \geq 1$ et $N \geq 0),$ on a
%
$$\card(E_{D,N})\leq  D\, 2^{D-1}\, e^{DN} \left(2^D\, e^{DN}+1 \right)^D.$$  

Montrons que 
$D\, 2^{D-1}\, e^{DN} \left(2^D\, e^{DN}+1 \right)^D \leq e^{D(D+1)(N+1)}$, 
ce qui donnera la majoration annonc\'ee. On a 
$$D\, 2^{D-1}\, e^{DN} \left(2^D\, e^{DN}+1 \right)^D = D\, 2^{D-1}\, e^{DN} \left( (2^D+1)e^{DN}-e^{DN}+1 \right) ^D$$
$$\leq D\, 2^{D-1}\, e^{DN} \left( (2^D+1)e^{DN} \right) ^D = D\, 2^{D-1}(2^D+1)^D\, e^{D(D+1)N}.$$
Il nous suffit donc de montrer que  $D\, 2^{D-1}(2^D+1)^D \leq e^{D(D+1)}$.
Par r\'ecurrence sur $D,$ on montre que pour tout $D \geq 2$, on a  
$D\,2^{D-1} \leq e^D $, 
et une \'etude facile de la fonction $x \mapsto e^x-2^x-1,$ montre que, pour tout $D \geq 2$,  
$(2^D+1)^D \leq e^{D^2}$. 
Pour $D=1$ l'in\'egalit\'e se montre par un calcul direct.

%
%

\subsubsection{Minoration de $\card(E_{D,N})$.}

Si $D \geq 1$ et $0 \leq N<1$, comme l'ensemble $E_{D,N}$ contient 0, alors 
$$e^{D(D+1)(N-1)} < 1 \leq \card(E_{D,N}).$$

Pour tout entier $D \geq 2$ et tout nombre r\'eel $N \geq 1,$ nous montrons que l'ensemble $E_{D,N}$ contient $\mathcal{A}_{D,H}$ pour $H=\frac{e^{DN}}{\sqrt{D+1}}.$
Soit $H'$ un nombre r\'eel positif, et soit $\alpha \in \mathcal{A}_{D,H'}.$ Alors $\deg(\alpha)= D$ et $\mathcal{H}(\alpha)\leq H'.$ D'apr\'es l'in\'egalit\'e (\ref{hauteurusuelle_hauteurlog}), on a 
$$h(\alpha)\leq \frac{1}{D}\log \mathcal{H}(\alpha)+\frac{1}{2D}\log(D+1)\leq \log \left( (D+1)^\frac{1}{2D}\, {H'}^\frac{1}{D} \right).$$
On en d\'eduit que $\mathcal{A}_{D,H'}\subset E_{D,N}$ d\`es que $\log ( (D+1)^{\frac{1}{2D}}H'^{\frac{1}{D}}) \leq N$. 

Comme $D \geq 2$ et $N \geq 1,$ alors $H=\frac{e^{DN}}{\sqrt{D+1}} \geq 2$ et nous pouvons appliquer la minoration obtenue pour le cardinal de $\mathcal{A}_{D,H}$ au lemme \ref{card A_{D,H}}: 
$$\card(E_{D,N})\geq \card(\mathcal{A}_{D,H}) > \frac{D}{8}(H-2)^{D+1}.$$
Or $\frac {e^D}{\sqrt{D+1}} \geq 4,$ donc $H \geq 4e^{D(N-1)} \geq 2 + 2e^{D(N-1)},$ d'o\`u $H-2 \geq 2e^{D(N-1)}$ et
$$\card(E_{D,N}) > \frac{D}{8} \left( 2e^{D(N-1)} \right)^{D+1} = D\,2^{D-2} e^{D(D+1)(N-1)}, $$
donc
$$\card(E_{D,N}) > e^{D(D+1)(N-1)}.$$

Si $\alpha$ est de degr\'e 1, on a $\alpha = a/b$ o\`u $a$ et $b \ne 0$ sont des entiers premiers entre eux. On a donc
$$h(\alpha)= \log(H(\alpha)) = \log( \max(|a|,|b|)) = \log(\mathcal{H}(\alpha))$$
et donc, pour tout $N \geq 1$, on a $E_{1,N} =\{ \alpha \in \mathbf{Q} ; \hspace{0,2cm} h(\alpha) \leq N\} = \{\alpha \in \mathbf{Q} ; \hspace{0,2cm} \log(\mathcal{H} (\alpha)) \leq N\}.$
On voit que $E_{1,N} \supset A_{1,H'}$ d\`es que $H' \leq e^N,$ donc, en particulier,
$$ \card(E_{1,N}) \geq \card(A_{1,H}) > (H-1)^2 \hspace{0,2cm} \mathrm{pour} \hspace{0,2cm} H=e^N.$$
Comme $ N\geq 1,$ alors $H = e^N \geq 2e^{N-1} \geq 1+e^{N-1},$ et on a 
$$ \card(E_{1,N}) > e^{2(N-1)},$$
ce qui compl\`ete la d\'emonstration du lemme \ref{e_{D,N}}.
\hfill $\Box$


\subsection{Autres r\'esultats.}\label{autres-resultats}

Dans ce travail, nous nous  servons du lemme \ref{e_{D,N}}, mais d'autres \'enonc\'es sont connus. Pour les pr\'esenter nous \'enon\c cons d'abord une cons\'equence du lemme \ref{e_{D,N}}.

Pour tout entier $d\geq 1$ et tout nombre r\'eel $N\geq 0,$ notons $\mathcal{E}_{d,N}$ l'ensemble des nombres alg\'ebriques de degr\'e {\it {exactement}} $d$ et de hauteur logarithmique absolue $\leq N$. En utilisant le lemme \ref{e_{D,N}}, on obtient l'encadrement suivant.

\begin{lemme}\label{mathcalE_d,N}

Pour tout entier $d \geq 2$ et tout nombre r\'eel $N \geq 0,$ le cardinal de $\mathcal{E}_{d,N}$ v\'erifie
$$c_1(d,N) e^{d(d+1)N} < \card(\mathcal{E}_{d,N}) < c_2(d,N) e^{d(d+1)N}$$
o\`u $c_1(d,N) = e^{-d(d+1)} - e^{-2dN + d^2 - d} \esp$ et $\esp c_2 (d,N) = (1 - e^{-2d(N+d)})e^{d(d+1)}.$

\end{lemme}

(Si $d=1$ et $N \geq 0,$ par un calcul simple analogue \`a celui de la fin de la preuve du lemme \ref{e_{D,N}}, on obtient $(e^N-1)^2 < \card(\mathcal{E}_{1,N}) \leq 2e^{2N} + e^N$.)

Remarquons que pour $N$ fixe et $d \to \infty$, la minoration du lemme \ref{mathcalE_d,N} est triviale car pour $d$ assez grand $c_1(d,N) < 0$. En revanche, T. Loher \cite{loher} (cf. aussi \cite{masser-vaaler}) a obtenu le r\'esultat suivant.

\begin{thm}{(Loher)}\label{Lo6}
Pour tout entier $d\geq 1$, et tout nombre r\'eel $N \geq \frac{1}{d}\log 2$,
$$\card(\mathcal{E}_{d,N}) \geq c_3(d) \,e^{d(d+1)N}$$
o\`u $c_3(d) = 2^{-d-2} (d+1)^{-\frac{1}{2}(d+1)}.$

\end{thm}

Pour $d$ fixe et $N \to \infty$, la minoration du lemme \ref{mathcalE_d,N} donne le m\^eme ordre de grandeur que le th\'eor\`eme de Loher, \`a savoir, $e^{d(d+1)N}$, mais la constante du th\'eor\`eme \ref{Lo6} est meilleure pour tout $d \geq 2.$ D'ailleurs, en remarquant que l'ensemble $E_{D,N}$ est la r\'eunion disjointe, pour $d \in \{1, \ldots , D\}$, des ensembles $\mathcal{E}_{d,N}$, on voit que le th\'eor\`eme de Loher implique: 

\textit{Pour tout entier $D\geq 1,$ et tout nombre r\'eel $N \geq \frac{1}{D}\log 2,$
$$\card(E_{D,N}) \geq c_3(D)\, e^{D(D+1)N},$$}
ce qui est meilleur que  la borne inf\'erieure du lemme \ref{e_{D,N}}, pour tout $D \geq 2.$ 

\bigskip

Pour la majoration, W. Schmidt obtient, dans \cite{schmidtI} (\'equation 1.4  p.170), une borne ayant m\^eme ordre de grandeur que celle du lemme \ref{mathcalE_d,N}, mais dont la constante est de qualit\'e l\'eg\`erement inf\'erieure  \`a $c_2(d,N)$. Il montre que {\it             
pour tout entier $d\geq 1,$ et tout nombre r\'eel $N \geq 0$,
$$\card(\mathcal{E}_{d,N}) \leq 2^{2d^2 + 14d + 11}\,e^{d(d+1)N}.$$
}
Concernant des estimations asymptotiques en $N$, il y a deux r\'esultats connus. Pour $d=1$, S.H. Schanuel (cf. \cite{schanuel} et aussi \cite{schmidtI} et \cite{masser-vaaler}) montre que
$$\card(\mathcal{E}_{1,N}) = \frac{12}{\pi^{2}}\, e^{2N} + O(Ne^{N}),$$
et pour $d=2$, W. Schmidt montre dans \cite{schmidtII}, que
$$\card(\mathcal{E}_{2,N}) = \frac{8}{\zeta(3)}\,e^{6N} + O(Ne^{4N}),$$
o\`u $\zeta$ est la fonction zeta de Riemann.
Aucune estimation asymptotique concernant le cardinal de $\mathcal{E}_{d,N}$ pour $d\geq 3$, ne semble actuellement connue (\cite{schmidt.dio} p.\,27 et \cite{schmidtII} p.\,346). N\'eanmoins, D.W. Masser et J. Vaaler \cite{masser-vaaler} ont obtenu le r\'esultat suivant.

\begin{thm}{(Masser-Vaaler)}\label{masser-vaaler}
Pour tout entier $d\geq 1$,  on a 
$$\lim_{N\rightarrow \infty} e^{-d(d+1)N}\card(\mathcal{E}_{d,N})  = c_4(d)\, $$
o\`u $c_4(d)=\frac{d\gamma(d)}{2\zeta(d+1)} \esp \textrm{et} \esp \gamma(d)=2^{d+1}(d+1)^{\delta} \prod_{k=1}^{\delta} \frac{(2k)^{d-2k}}{(2k+1)^{d+1-2k}}$ avec $\delta=\left[ \frac{d-1}{2} \right].$

\end{thm}
(La fonction $\gamma(d)$ v\'erifie $\lim_{d \rightarrow \infty} \frac{\log \gamma(d)}{d \log d} = -\frac{1}{2}$.)
Le th\'eor\`eme \ref{masser-vaaler} implique le r\'esultat suivant.
{\textit Pour tout entier $D\geq 1$, il existe un entier $N_0$ tel que, pour tout $N \geq N_0$, on a 
$$\card(E_{D,N}) \leq \frac{2^{D-1}\,D^2\, (D+1)^{\frac{D-1}{2}}}{\zeta(D+1)}\, e^{D(D+1)N}.$$
}


\section{Construction d'exemples}\label{exemples}

Dans cette section, nous construisons la fonction $f$ du th\'eor\`eme \ref{contrex}. 
Nous construisons aussi (th\'eor\`eme \ref{contrex2}) une fonction $g$ enti\`ere et transcendante, dont toutes ses d\'eriv\'ees, envoient tout nombre alg\'ebrique $\alpha$ dans $\mathbf{Z}\left[\frac{1}{2},\alpha\right]$.


\subsection{Construction de la fonction $f$.}\label{exemplef}

On rappelle que $\epsilon_{D,N}= \card \{ \alpha \in \Qbarre / \espa [\rat (\alpha):\rat ] \leq D , \esp h(\alpha )\leq N \}$.

Soit $\phi$ une fonction positive telle que $\phi(x)/x$ tende vers $0$ quand $x$ tende vers l'infini.  

Nous nous donnons une suite $(b_k)_{k\geq 1}$ de nombres r\'eels $>0$ telle que la s\'erie $\sum_{k\geq1}b_k$ soit convergente, et un nombre r\'eel $x_0 \geq 1,$ tel que pour tout nombre r\'eel $x \geq x_0,$
\begin{equation}
\phi(x) \leq x-1. \label{x_0}
\end{equation}

Nous allons construire, r\'ecursivement, une suite strictement croissante  $(N_{\delta})_{\delta\geq 1}$ de nombres r\'eels $\geq x_0$ tendant vers l'infini et une suite $(a_{\delta})_{\delta \geq 1}$ de nombres rationnels v\'erifiant les conditions suivantes.

\bigskip

\noindent i)  Pour une infinit\'e de $\delta,$ on a $a_{\delta} \ne 0$ et pour tout $\delta \geq 1,$
\begin{equation} 
|a_{\delta}| \leq b_{\delta}\,(\delta + e^{\delta N_{\delta}})^{-\delta \epsilon_{\delta,N_{\delta}}}.  \label{ak}
\end{equation}

\noindent ii) Pour tout $\delta \geq 2,$ 
\begin{equation}
N_{\delta} \geq 2 \left( \log(\delta-1) + \sum_{k=1}^{\delta-1} h(a_k) + (\delta-1)^2 \,\epsilon_{\delta-1,N_{\delta-1}}(\log2 + 1 + N_{\delta-1}) \right).  \label{Ndelta>.} 
\end{equation}

\noindent iii) Pour tout $\delta \geq 2,$
\begin{equation}
\frac{N_{\delta}}{\phi(N_{\delta})} \geq 2 (\delta-1)^2 \,\epsilon_{\delta-1,N_{\delta-1}}.                                                   \label{Ndelta>2 Phi}
\end{equation}

\bigskip


Pour la suite $(b_k)_{k\geq1}$ on peut prendre, par exemple, $b_k = 2^{-k}$, pour tout $k \geq 1.$

Pour construire $(a_{\delta}, N_{\delta})_{\delta \geq 1}$, on proc\`ede de la fa\c con suivante.

On pose $N_1 = [x_0] + 1, \hspace{0,2cm} c_1 =  1 + \left[ 2 (1 + e^{N_1})^{\epsilon_{1,N_1}} \right],$ et $a_1 = \frac{1}{c_1}.$ 

Soit $\delta \geq 2.$ Pour $1 \leq k \leq \delta-1,$ on suppose d\'efinis $a_k$ et $N_k$ v\'erifiant, pour $k \in \{1, \ldots, \delta - 1\},$
$$|a_k| \leq b_k\,(k + e^{k N_k})^{-k \,\epsilon_{k,N_k}},$$
et pour $k \in \{ 2, \ldots, \delta - 1\} ,$
$$N_k \geq 2 \left( \log(k-1) + \sum_{r=1}^{k-1} h(a_r) + (k-1)^2 \,\epsilon_{k-1,N_{k-1}}(\log2 + 1 + N_{k-1}) \right)$$
et
$$\frac{N_k}{\phi(N_k)} \geq 2 (k-1)^2 \,\epsilon_{k-1,N_{k-1}}.$$

(Remarquer que la derni\`ere in\'egalit\'e est possible puisque la fonction $\phi(x)/x$ tend vers 0 quand $x$ tend vers l'infini.)

On choisit pour $N_{\delta}$ un nombre r\'eel v\'erifiant les conditions (\ref{Ndelta>.}) et (\ref{Ndelta>2 Phi}), et on pose $c_{\delta} = 1 + \left[ 2^{\delta} (\delta + e^{\delta N_{\delta}})^{\delta \,\epsilon_{\delta, N_{\delta}}} \right]$ et  $a_{\delta} = c_{\delta}^{-1}.$ 
 
Il est clair que les $a_{\delta}$ peuvent \^etre choisis de fa\c con \`a ce que tous (\`a partir d'un ceratin rang) soient non nuls. 

Ainsi, la suite $(a_{\delta}, N_{\delta})_{\delta\geq1}$ satisfait les conditions (\ref{ak}), (\ref{Ndelta>.}) et (\ref{Ndelta>2 Phi}). En particulier, la condition (\ref{Ndelta>.}) implique que la suite $(N_{\delta})_{\delta\geq1}$ est strictement croissante et tend vers l'infini.


On pose, pour tout $k \geq 1,$ 
$$P_k(X) = \prod_{\beta \in E_{k,N_k}}(X - \beta )^k.$$

On d\'efinit la fonction $f$ en posant, pour tout $z\in \complexes,$ 
$$f(z) = \sum_{k\geq 1}a_k P_k(z).$$

\bigskip

Pour d\'emontrer le th\'eor\`eme \ref{contrex}, nous aurons besoin des lemmes suivants.


\begin{lemme}\label{entiere-et-transc}

La fonction $f$ est une fonction enti\`ere et transcendante sur $\complexes (z).$

\end{lemme}

\noindent \textbf{D\'emonstration du lemme \ref{entiere-et-transc}.}

Soient $R>0$ et $z\in \DRferme.$ Nous \'ecrivons la s\'erie   
$$\sum_{k\geq1} a_k P_k(z) = \sum_{k=1}^{[R]} a_k P_k(z) + \sum_{k\geq [R] + 1} a_k P_k(z).$$
La premi\`ere somme est finie et d\'efinit un polyn\^ome. Nous regardons donc la deuxi\`eme. Pour tout $k\geq [R]+1,$ nous avons
$$ |a_k P_k(z)| = |a_k| \prod_{\alpha \in E_{k,N_k}} |z - \alpha|^k \leq |a_k| \prod_{\alpha \in E_{k,N_k}} (|z| + |\alpha|)^k.$$
Pour $\alpha \in E_{k,N_k},$ nous avons $|\alpha| \leq M(\alpha) = e^{\deg(\alpha) h(\alpha)} \leq e^{k N_k},$
et donc
$$|a_k P_k(z)|  \leq  |a_k| \prod_{\alpha \in E_{k,N_k}} (|z| + e^{k N_k})^k.$$
Or $|z| \leq R \leq k$, donc $|a_k P_k(z)| \leq |a_k| \left( k + e^{k N_k}\right)^{k \epsilon_{k, N_k}}$, 
et, d'apr\`es la condition (\ref{ak}) sur la suite $(a_{\delta})_{\delta \geq 1}$, 
$$|a_k P_k(z)| \leq b_k.$$

Comme la s\'erie $\sum b_k$ est convergente, ceci montre que la s\'erie $\sum_{k\geq1}a_k P_k(z)$ converge normalement sur tout compact de $\complexes$; sa limite $f$ est une fonction holomorphe sur $\complexes$ tout entier.

%
%

Montrons maintenant que la fonction $f$ n'est pas polynomiale; comme elle est enti\`ere, cela montrera qu'elle est transcendante.

Supposons que $f$ est un polyn\^ome de degr\'e $n \geq 1$. 

Comme la suite $(\epsilon_{k,N_k})_{k \geq 1}$ est strictement croissante, il existe $k_0 \geq 1$, tel que pour tout $k \geq k_0$, on a $(k-1) \epsilon_{k-1,N_{k-1}} \geq n$.

Soit $K \geq k_0$ tel que $a_{K-1}\ne 0$. Consid\'erons la somme 
$$g_K(z) = \sum_{k=1}^{K-1} a_k P_k(z).$$
C'est un polyn\^ome de degr\'e $(K-1) \epsilon_{K-1,N_{K-1}} \geq n$. La diff\'erence
$$f(z) - g_K(z) = \sum_{k\geq K} a_k P_k(z)$$
est un polyn\^ome de degr\'e $\leq \max\{(K-1) \epsilon_{K-1,N_{K-1}}, n  \} = (K-1) \epsilon_{K-1,N_{K-1}}$. Or, pour tout $k\geq K,$ le polyn\^ome $P_k$ est multiple de $P_K,$ donc $f(z) - g_K(z)$ s'annule en tous les z\'eros de $P_K,$ c'est-\`a-dire, en $K \epsilon_{K,N_K}$ points.

Comme le polyn\^ome $f(z) - g_K(z)$ est de degr\'e $\le (K-1) \epsilon_{K-1,N_{K-1}} <K \epsilon_{K,N_K}$, c'est le polyn\^ome nul, et donc
$$-a_K P_K(z) = \sum_{k \geq K+1} a_k P_k(z).$$
Soit $z_0 \in E_{K+1,N_{K+1}} \setminus E_{K,N_K}$. Alors 
$$\forall k \geq K+1, \esp P_k(z_0) = 0 \esp \textrm{et} \esp  P_K(z_0) \ne 0.$$
Nous en concluons que, pour tout $K \geq k_0$, on a $a_K=0,$ ce qui contredit la condition i).\hfill $\Box$



\begin{lemme}\label{f(alpha)}

Pour tout entier naturel $\sigma,$ la fonction d\'eriv\'ee $f^{(\sigma)}$ envoie tout nombre alg\'ebrique $\alpha,$ dans $\rat(\alpha).$

%

\end{lemme}

\noindent \textbf{D\'emonstration du lemme \ref{f(alpha)}.}

Soient $\alpha \in \Qbarre$ et $\sigma \in \naturels.$
Par hypoth\`ese, la suite $(N_{\delta})_{\delta \geq 1}$ tend vers l'infini, donc
$$\Qbarre = \bigcup_{k\geq1} E_{k,N_k},$$
et comme les ensembles $(E_{k,N_k})_{k\geq 1}$ forment une suite croissante pour l'inclusion, il existe $k' \geq 1$ tel que $\alpha \in E_{k,N_k}$, pour tout $k \geq k'$. 
Notons $k_0 = \min\{k'\geq 1/ \esp \forall k \geq k', \esp \alpha \in E_{k,N_k} \}$ et $M=\max\{ k_0, \sigma+1 \}.$ 

Si $M=1$ (i.e. $\sigma = 0$ et $k_0 = 1$), alors $f(\alpha) = 0$.

Si $M>1$, pour tout $k \geq M,$ on a $\alpha \in E_{k, N_k}$ et $k>\sigma.$ 
%
%
Donc, pour tout $k \geq M, \esp P_k^{(\sigma)}(\alpha) = 0$ et 
$$f^{(\sigma)}(\alpha) = \sum_{k=1}^{M - 1}a_k P_k^{(\sigma)}(\alpha).$$

De plus, les polyn\^omes $P_k$ sont des puissances de polyn\^omes unitaires dont les ensembles de z\'eros sont des r\'eunions de syst\`emes complets de conjugu\'es sur $\bf{Q}$, ils sont donc fix\'es par tout \'el\'ement du groupe de Galois de $\Qbarre$ sur $\rat.$ Ils sont donc \`a coefficients rationnels, ainsi que tous leurs polyn\^omes d\'eriv\'es, et comme, par d\'efinition, les nombres $a_k$ sont aussi rationnels,
$$f^{(\sigma)}(\alpha) \in \rat(\alpha).$$
\hfill $\Box$


Le lemme suivant sera lui aussi utile pour la d\'emonstration du th\'eor\`eme \ref{contrex}.


\begin{lemme}\label{un-demi}

Pour tout entier $D \geq 1$ et tout nombre r\'eel $N \geq 0,$ nous avons
$$\card \left(E_{D,N} \cap \Duferme \right) \geq \frac{1}{2} \card(E_{D,N}).$$

\end{lemme}

\noindent \textbf{D\'emonstration du lemme \ref{un-demi}.}

Soient $D$ un entier positif et $N$ un nombre r\'eel positif ou nul.

Pour un nombre alg\'ebrique $\alpha$ non nul, on a 
$\deg(\alpha) = \deg \left(\frac{1}{\alpha}\right)$ et $h(\alpha) = h\left(\frac{1}{\alpha}\right)$,
donc $\alpha$ appartient \`a $E_{D,N}$ si et seulement si $\frac{1}{\alpha}$ appartient \`a $E_{D,N}$. De plus, $\alpha$ appartient au disque ferm\'e $\Duferme$ si et seulement si $\frac{1}{\alpha}$ n'appartient pas au disque ouvert $\Duouvert$.
D'o\`u 
$$\card(E_{D,N} \cap \Duferme ) = \card(E_{D,N} \setminus \Duouvert ) + 1,$$
et donc
\begin{displaymath}
\begin{array}{rcl}
\card(E_{D,N}) & =  & \card(E_{D,N} \cap \Duferme ) + \card(E_{D,N} \setminus \Duouvert ) \\
               &  \leq   & 2\, \card(E_{D,N} \cap \Duferme ),
\end{array}
\end{displaymath}
ce qui d\'emontre le lemme. \hfill $\Box$

\bigskip

\noindent \textbf{D\'emonstration du th\'eor\`eme \ref{contrex}.}

Soient $D$ et $d$ des entiers tels que $d \geq D  \geq 1$.

Nous commen\c cons par remarquer que, d'apr\`es le lemme \ref{f(alpha)}, 
$$\Sigma_{D,N_d}(f,1)= \{\alpha \in \Qbarre \cap \Duferme / \esp \deg(\alpha) \leq D, \esp h(\alpha) \leq N_d, \esp h(f(\alpha)) \leq N_d \}.$$

Montrons que cet ensemble contient $E_{D,\phi(N_d)+1} \cap \Duferme.$
Soit $\alpha \in \Qbarre$ tel que 
$$\deg(\alpha) \leq d \esp \textrm{et} \esp  h(\alpha) \leq \phi(N_d) + 1.$$

Comme $N_d \geq x_0,$ alors, d'apr\`es (\ref{x_0}), $\phi(N_d) \leq N_d - 1,$ et $h(\alpha) \leq N_d.$
Comme, en plus, $\deg (\alpha) \leq d,$ alors $\alpha \in E_{d, N_d}$ et pour tout $k \geq d$, on a $P_k(\alpha) = 0.$

Si $d=D=1,$ alors $f(\alpha) =0$ et on a l'inclusion $E_{1,\phi(N_1)+1} \cap \Duferme \subset \Sigma_{1,N_1}(f,1).$  

Supposons maintenant que $d \ge D \geq 1$ et $d\ge 2.$ Alors
$$f(\alpha) = \sum_{k=1}^{d-1} a_k P_k(\alpha).$$
Majorons sa hauteur. En appliquant les formules (\ref{h(a_1+...+a_d)<h(a_1)+...+h(a_d)}), puis (\ref{h(a.b)<h(a)+h(b)}), on obtient 
$$h(f(\alpha))  =  h\left(\sum_{k=1}^{d-1} a_k P_k(\alpha)\right) \leq \log(d-1) + \sum_{k=1}^{d-1} h(a_k P_k(\alpha))$$
$$\leq  \log(d-1) + \sum_{k=1}^{d-1} h(a_k) + \sum_{k=1}^{d-1} h(P_k(\alpha)).$$
Soit $k \in \{1, \ldots, d-1\}.$ Nous avons
$$P_k(\alpha) = \prod_{\beta \in E_{k, N_k}} (\alpha - \beta)^k.$$
En utilisant les formules (\ref{h(a.b)<h(a)+h(b)}) et (\ref{h(a_1+...+a_d)<h(a_1)+...+h(a_d)}), on obtient 
$$h(P_k(\alpha)) \leq k\sum_{\beta \in E_{k, N_k}} h(\alpha - \beta) \leq k\sum_{\beta \in E_{k, N_k}} (\log2 + h(\alpha) + h(\beta)),$$
et comme, pour tout $\beta \in E_{k, N_k}$, on a $h(\beta) \leq N_k$, et que, par hypoth\`ese, $h(\alpha) \leq \phi(N_d) + 1$, alors
$$h(P_k(\alpha)) \leq k\,\epsilon_{k, N_k} ( \log 2 + \phi(N_d) + 1 + N_k ).$$
D'o\`u
$$h(f(\alpha)) \leq  \log(d-1) + \sum_{k=1}^{d-1} h(a_k) + \sum_{k=1}^{d-1} k\,\epsilon_{k, N_k} ( \log 2 + \phi(N_d) + 1 + N_k ).$$
Or, les suites $(N_k)_{k\geq1}$ et $(\epsilon_{k, N_k})_{k\geq1}$ sont croissantes, donc $h(f(\alpha))$ est inf\'erieure \`a
$$\log(d-1) + \sum_{k=1}^{d-1} h(a_k) + (d-1)^2 \epsilon_{d-1, N_{d-1}}( \log 2 + \phi(N_d) + 1 + N_{d-1})  $$   
$$=\phi(N_d)\left( (d-1)^2 \epsilon_{d-1, N_{d-1}} \right) + \log(d-1) + \sum_{k=1}^{d-1} h(a_k) + (d-1)^2 \epsilon_{d-1, N_{d-1}} (\log 2 + 1 + N_{d-1})$$ 
ce qui est $\leq N_d$, d'apr\`es les conditions (\ref{Ndelta>.}) et (\ref{Ndelta>2 Phi}).

En particulier, si $\alpha$ appartient \`a $E_{D,\phi(N_d)+1} \cap \Duferme$, alors $\deg(\alpha) \leq D \leq d$ et $h(f(\alpha)) \leq N_d$, ce qui montre que $\alpha$ appartient \`a $\Sigma _{D, N_d}(f,1)$. 
Ainsi, pour tout couple $(D, d)$ avec $d \geq D \geq 1$, nous avons
$$\sigma_{D,N_d} = \card(\Sigma_{D,N_d} (f,1)) \geq \card(E_{D,\phi(N_d)+1} \cap \Duferme).$$
En appliquant le lemme \ref{un-demi} et le lemme \ref{e_{D,N}}, nous obtenons
$$\sigma_{D,N_d} \geq \frac{1}{2}\, \card(E_{D,\phi(N_d)+1}) > \frac{1}{2} \,e^{D(D+1)\phi(N_d)}.$$ \hfill $\Box$

\bigskip

On remarque que la d\'emonstration se simplifie consid\'erablement si on se contente de construire $f$ enti\`ere et transcendante v\'erifiant (\ref{a)}). Pour cela, il suffit de se donner une suite strictement croissante  $(N_{\delta})_{\delta\geq 1}$ de nombres r\'eels $\geq 1$ tendant vers l'infini et une suite $(a_{\delta})_{\delta \geq 1}$ de nombres rationnels v\'erifiant la condition i).

Le th\'eor\`eme \ref{contrex}, qui fait intervenir les d\'eriv\'ees de la fonction transcendante, pose le probl\`eme suivant. Peut-on esp\'erer avoir un \'enonc\'e dans le sens du th\'eor\`eme \ref{propprincipale}, faisant lui aussi, intervenir les d\'eriv\'ees?

%
%
%

\subsection{Construction de la fonction $g$.}\label{exempleg}

En modifiant la construction de la fonction $f$ du th\'eor\`eme \ref{contrex}, on peut construire une fonction $g$ enti\`ere et transcendante, v\'erifiant la m\^eme conclusion (\ref{b)}) concernant le cardinal de $\Sigma_{D,N}(g,1)$, mais pour laquelle la premi\`ere conclusion est chang\'ee. Pr\'ecis\'ement :

\begin{thm}\label{contrex2}

Il existe une fonction $g$ enti\`ere et transcendante sur $\mathbf{C}(z)$ telle que,
\begin{equation}
\forall \alpha \in \Qbarre, \esp \forall \sigma \geq 0, \esp g^{(\sigma)}(\alpha) \in \mathbf{Z}\left[\frac{1}{2},\alpha \right]. \label{a bis)}
\end{equation}

\end{thm}

\noindent{\bf D\'emonstration du th\'eor\`eme \ref{contrex2}.}

On se donne une suite croissante $(N_k)_{k\geq1}$ de nombres r\'eels $\geq 1$ tendant vers l'infini et une suite $(b_k)_{k\geq1}$ de nombres rationnels tels que :

\bigskip

\noindent i) pour une infinit\'e de $k\geq 1, \esp b_k \ne 0$,

\medskip

\noindent ii) pour tout $k\geq1$, on ait $|b_k| \leq 2^{-k}(k+e^{k N_k})^{-k \epsilon_{k,N_k}}\ $, 

\medskip

\noindent iii) pour tout $k\geq1$,
$b_k \in \mathbf{Z} \left[\frac{1}{2}\right]$. 
%
(Comme $\enteros\left[\frac{1}{2}\right]$ est dense dans $\rat$, ceci est possible.)
\medskip

On pose, pour tout $k\geq 1$, 
$$Q_k(X) = \prod_{\beta \in E_{k,N_k}}P_{\beta}(X)^k,$$
o\`u $P_{\beta}$ est le polyn\^ome minimal de $\beta$ sur $\enteros$,
et, pour tout $z \in \mathbf{C}$,
$$g(z) = \sum_{k\geq1}b_k Q_k(z).$$

\bigskip

Les lemmes suivants d\'emontrent le th\'eor\`eme \ref{contrex2}.

\begin{lemme}
La fonction $g$ est enti\`ere et transcendante.
\end{lemme}

La d\'emonstration est identique \`a celle du lemme \ref{entiere-et-transc}.

\begin{lemme}\label{g(alpha)}
Pour tout entier naturel $\sigma,$ la fonction d\'eriv\'ee $g^{(\sigma)}$ envoie tout nombre alg\'ebrique $\alpha$ dans $\mathbf{Z}\left[\frac{1}{2},\alpha\right].$

\end{lemme}

\noindent \textbf{D\'emonstration du lemme \ref{g(alpha)}.}

Cette d\'emonstration suit celle du lemme \ref{f(alpha)}.
On se donne  $\alpha$ dans $\Qbarre$ et $\sigma$ dans $\mathbf{N}$.
Il existe un nombre entier $M\geq 1$, tel que, pour tout $k\geq M$, on ait $\alpha \in E_{k,N_k}$ et $k>\sigma$, donc $Q_k^{(\sigma)} (\alpha) = 0$ et
$$g^{(\sigma)}(\alpha) = \sum_{k=1}^{M - 1}b_k Q_k^{(\sigma)}(\alpha).$$

Comme les polyn\^omes $Q_k$ sont \`a coefficients dans $\mathbf{Z}$ et que, de plus, d'apr\`es la condition iii), pour tout $k\geq 1,\esp b_k \in \mathbf{Z}\left[ \frac{1}{2} \right]$, on a
$g^{(\sigma)}(\alpha) \in \mathbf{Z}\left[ \frac{1}{2},\alpha \right]$,
ce qui d\'emontre le lemme \ref{g(alpha)} et donc le th\'eor\`eme \ref{contrex2}. 
\hfill $\Box$

\bigskip

Si on se donne une fonction $\phi$ comme au th\'eor\`eme \ref{contrex} et si les suites $(N_k)$ et $(b_k)$ intervenant dans la fonction $g$ v\'erifient de plus les conditions (\ref{Ndelta>.}) (avec $(b_k)$ \`a la place de $(a_k)$) et (\ref{Ndelta>2 Phi}), alors la fonction $g$ ainsi obtenue satisfait de plus la conclusion (\ref{b)}) du th\'eor\`eme \ref{contrex}.

\end{document}